\documentclass{amsart}

\usepackage{amscd,amssymb,amsmath,graphicx,verbatim}
\usepackage[dvips]{hyperref}
\usepackage[TS1,OT1,T1]{fontenc}
\usepackage{xypic}
\newtheorem{theorem}{Theorem}[section]
\newtheorem{lemma}[theorem]{Lemma}
\newtheorem{corollary}[theorem]{Corollary}
\newtheorem{proposition}[theorem]{Proposition}

\theoremstyle{definition}
\newtheorem{definition}[theorem]{Definition}
\newtheorem{example}[theorem]{Example}

\theoremstyle{remark}
\newtheorem{remark}[theorem]{Remark}

\numberwithin{equation}{section}

\begin{document}
\title{Mutations vs. Seiberg Duality} 
\author{Jorge Vitória}
\maketitle
\begin{abstract}
For a quiver with potential, Derksen, Weyman and Zelevinsky defined a combinatorial transformation - mutations. Mukhopadhyay and Ray, on the other hand, tell us how to compute Seiberg dual quivers for some quivers with potentials through a tilting procedure, thus obtaining derived equivalent algebras. 
In this text, we compare mutations with the concept of Seiberg duality given by [10], concluding that for a certain class of potentials (the good ones) mutations coincide with Seiberg duality, therefore giving derived equivalences.
\end{abstract}
\begin{section}{Preliminaires}
In this section we introduce the material from [7] that will be used and recall some definitions.

We will use the following notation: $\mathbb{K}$ is a field; $\mathbb{K}Q$ is the path algebra of the quiver $Q$ over $\mathbb{K}$ (concatenation of paths is written as composition of functions); $Proj(R)$ is the full subcategory of projective right modules over a $\mathbb{K}$-algebra $R$; $P(R)$ is the full subcategory of finitely generated projective right modules over $R$; $K^{b}(Q)$ and $D^{b}(Q)$ are, respectively, the bounded homotopy category and the bounded derived category of right modules over  $\mathbb{K}Q$.\\

\begin{definition}
A \textbf{potential} on a quiver is an element of the vector space spanned by the cycles of the quiver (denote it by $\mathbb{K}Q_{cyc}$).
\end{definition}
\begin{definition}Let $A = <Q_{1}>$, i.e., the vector space spanned by all arrows. For each $\xi \in A^{*}$ (the dual of $A$), define a \textbf{cyclic derivative}: 
	\begin{equation}
	\partial/\partial\xi:\begin{array}{rcl} \mathbb{K}Q_{cyc} & \to & \mathbb{K}Q \\ a_1 \ldots a_n  & \mapsto & \sum_{k=1}^{n} {\xi(a_{k})a_{k+1} \ldots a_{n}a_{1} \ldots a_{k-1}} \end{array}.
	\end{equation}
\end{definition}

If $x\in Q_{1}$, we will denote by $\partial/\partial x$ the cyclic derivative correspondent to the element of $A^{*}$ which is the dual of $x$ in the dual basis of $A$. \\

\begin{definition} 
Two potentials are  \textbf{cyclically equivalent} if $S - S'$ lies in the span of elements of the form $a_{1} \ldots a_{n-1}a_{n} - a_{2} \ldots a_{n}a_{1}$. A pair $(Q,S)$ is said to be a quiver with potential if $Q$ has no loops and no two cyclically equivalent paths appear on $S$. Two quivers with potentials $(Q,S)$ and $(\tilde{Q},\tilde{S})$ are said to be \textbf{right equivalent} if there is isomorphism $\phi$ between $\mathbb{K}Q$ and $\mathbb{K}\tilde{Q}$ such that $\phi(S)$ is cyclically equivalent to $\tilde{S}$.\\
\end{definition}

\begin{definition}
Given a quiver with potential $(Q,S)$, define the \textbf{jacobian algebra of $(Q,S)$} as $J(Q,S) = \mathbb{K}Q / <J(S)>$, where $J(S) = (\partial S/\partial x)_{x\in Q_{1}}$.\\
\end{definition}
\begin{remark}
Two right equivalent quivers with potentials have isomoprhic jacobian algebras (see [7]).\\
\end{remark}
\begin{definition} 
Define the \textbf{trivial part} of a quiver with potential $(Q_{triv},S_{triv})$ by taking $S_{triv}$ as the degree two homogeneous component of S and $Q_{triv}$ as the subquiver of $Q$ consisting only in the arrows appearing in $S_{triv}$. The \textbf{reduced part} $(Q_{red},S_{red})$ is formed by the non-trivial part of the potential $S$ and by the quiver obtained by taking the quotient of $A$ by the arrows appearing in $S_{triv}$.
\end{definition}

The following theorem will allow us to define mutation on a quiver with potential.\\

\begin{theorem}[7]
For a quiver with potential $(Q,S)$, there exist a trivial quiver with potential $(Q_{triv},S_{triv})$ and a reduced quiver with potential $(Q_{red},S_{red})$ such that $(Q,S)$ is right equivalent to $(Q_{triv}\oplus Q_{red},S_{triv} + S_{red})$ - $Q_{triv}\oplus Q_{red}$ stands for the quiver obtained by taking the direct sum of the arrow spans.
\end{theorem}

Let's now describe the procedure of mutation of a quiver with potential $(Q,S)$ on a vertex $k$ (denote it by $\mu_k(Q,S)$).

\begin{enumerate}
	\item Suppose $k$ does not belong to any 2-cycle and that $S$ doesn't have any cycle starting and finishing on $k$ (if it does, substitute it by a cyclically equivalent potential that doesn't).
	\item Change the quiver in the following way:
	\begin{itemize}
	\item Reflect arrows starting or ending at $k$. Denote reflected arrows by $(.)^{*}$;
	\item Create one new arrow for each path of the form \xymatrix{{\bullet}^{i} \ar[r]^{\alpha}&{\bullet}^{k}\ar[r]^{\beta}&{\bullet}^{j}} and denote it by $[\beta\alpha]$. We denote the resulting quiver by $\tilde{Q}$.
	\end{itemize}
\item Change the potential in the following way:
\begin{itemize}
	\item Substitute factors appearing in $S$ of the form $\beta\alpha$ by the new arrow $[\beta\alpha]$ and denote it by $[S]$;
	\item Add $\Delta_{k} = \sum_{\xymatrix{{\bullet}^{i} \ar[r]^{\alpha}&{\bullet}^{k}\ar[r]^{\beta}&{\bullet}^{j}}} {[\beta\alpha]\alpha^{*}\beta^{*}}$ to $[S]$. We denote the resulting potential by $\tilde{S}$.
\end{itemize}
\item The mutation at $k$ of $(Q,S)$ is $\mu_k(Q,S) =(\bar{Q},\bar{S}) := (\tilde{Q}_{red},\tilde{S}_{red})$
\end{enumerate}

Note that these mutations generalize reflection functors on quivers with no relations in the sense that if you do a mutation on either a sink or a source, this procedure reduces to reflect arrows on that vertex.

Let us recall Rickard's theorem, starting by defining tilting complex ([12]).\\

\begin{definition} A \textbf{tilting complex} over a ring $R$ is an object $T$ of $K^{b}(P(R))$, such that:
\begin{enumerate}
	\item $\forall i \neq 0$, $Hom_{K^{b}(P(R))}(T,T[i]) = 0$;
	\item T generates $K^{b}(P(R))$ as a triangulated category.\\
\end{enumerate}
\end{definition}

\begin{theorem}[Rickard]
Let $R$ and $S$ be two rings. Then $D^{b}(R)$ is equivalent to $D^b(S)$ iff there is a tilting complex $T$ over $R$ such that $S \cong End_{K^{b}(R)}(T)^{op}$.
\end{theorem}
\end{section}

\begin{section}{Seiberg Duality}
We will define Seiberg duality on quivers as a tilting procedure and therefore as an equivalence of derived categories. To check if a complex is tilting we will have to compute homomorphisms in the derived category between (finitely generated) projective modules. \\

\begin{remark}
Note that $K^{b}(P(R))$ is a full subcategory of $K^{b}(R)$ and therefore, for an object $T$ in $K^{b}(P(R))$ (in particular for a tilting complex) we have $End_{K^{b}(R)}(T)^{op} = End_{K^{b}(P(R))}(T)^{op}$.
\end{remark}

Let $(Q,S)$ be a quiver with potential with $n$ vertices such that every vertex is contained in some cycle (which we shall assume from now on) and, for each vertex $k$, consider the following complex:
\begin{equation}\nonumber
T^{k} = \oplus_{i=1}^{n} T^{k}_{i}
\end{equation}
where
\begin{equation}\nonumber
T^{k}_{i} = 0 \rightarrow P_{i} \rightarrow 0,\  if\  i \neq k
\end{equation}
and
\begin{equation}\nonumber
T^{k}_{k} = \xymatrix{{0} \ar[r] & {\oplus_{j \rightarrow k} P_{j}}  \ar[r]^{(\alpha_{j})_{j}} & {P_{k}} \ar[r]& {0}}
\end{equation}\\

\begin{lemma}
$T^{k}$ is a tilting complex over the jacobian algebra of $(Q,S)$ if and only if $Hom_{K(P(Q))}(T^{k}_{k},T^{k}_{s}[-1])=0$, $\forall s$.
\end{lemma}

\begin{proof}
\begin{enumerate}
\item $Hom_{K^{b}(P(Q))}(T^{k},T^{k}[i]) = 0$ $\forall i \neq 0$.
It is clear that if $r, s \neq k$, then $Hom_{K(P(Q))}(T^{k}_{r},T^{k}_{s}[i]) = 0, \forall i \neq 0$.
Now, suppose $s = k$ and $r \neq k$ . Then we only have to check that the set $Hom_{K(P(Q))}(T^{k}_{r},T^{k}_{k}[1])$ reduces to zero. Note that, since a homomorphism between $P_{r}$ to $P_{k}$ is identified with an element of the path algebra with each term being a path from $r$ to $k$, every such homomorphism factors through $\oplus_{j \rightarrow k} P_{j}$.
\begin{equation}\nonumber
\xymatrix{& {0}\ar[r] & {P_{i}}\ar[r]\ar@{.>}[ld]\ar[d] &{0} \\ {0}\ar[r] &{\oplus_{j \rightarrow k} P_{j}} \ar[r]& {P_{k}} \ar[r]& 0}
\end{equation}
Such factorization implies that these maps of complexes are homotopic to zero, thus zero in the homotopy category.

If $r = k$ then we also have such a homotopy just by taking identity maps.
\begin{equation}\nonumber
\xymatrix{& {0}\ar[r] & {\oplus_{j \rightarrow k} P_{j}}\ar[r]\ar@{.>}[ld]\ar[d] &{P_{k}}\ar[r] & {0} \\ {0}\ar[r] &{\oplus_{j \rightarrow k} P_{j}} \ar[r]& {P_{k}} \ar[r]& 0}
\end{equation}

\item add($T^{k}$) generates $K^{b}(P(Q))$ as a triangulated category.
It is enough to prove that the stalk complexes of indecomposable projective modules are generated by the direct summands of $T^{k}$. 

Consider the direct summands of $T^{k}$ and take the cone of the map $T^{k}_{k}$ to $\oplus_{j \rightarrow k} T^{k}_{j}$ defined by:
\begin{equation}\nonumber\xymatrix{{0}\ar[r] &{\oplus_{j\rightarrow k} P_{j}}\ar[r]\ar[d]_{id} & {P_{k}}\ar[r] & {0} \\ {0}\ar[r] & {\oplus_{j\rightarrow k} P_{j}}\ar[r] & {0}} 
\end{equation}
That cone is just the following complex (the underlined term is in degree zero):
\begin{equation}
\xymatrix{{0}\ar[r]& {\oplus_{j\rightarrow k} P_{j}} \ar[rr]^{((\alpha_{j})_{j},0)}&&{\underline{P_{k}\oplus (\oplus_{j\rightarrow k} P_{j})}}\ar[r]&{0} }
\end{equation}
Consider the map from the complex (2.1) to the stalk complex of $P_{k}$ in degree zero defined by identity in the first component and $-(\alpha_{j})_{j}$ in the second component and consider the map from this same stalk complex to (2.1) defined by the inclusion of $P_{k}$. We will prove that the composition of these maps is homotopic to the identity map, hence proving that these complexes are isomorphic in the derived category. In fact, that follows from the following diagram:
\begin{equation}\nonumber
\xymatrix{{0}\ar[rr]&& {\oplus_{j\rightarrow k} P_{j}}\ar@{.>}[ddddll]\ar[dd]\ar[rr]^{((\alpha_{j})_{j},0)}&&{\underline{P_{k}\oplus (\oplus_{j\rightarrow k} P_{j})}}\ar@{.>}[ddddll]^{(0,id)}\ar[dd]^{(id,-(\alpha_{j})_{j})}\ar[rr]&&{0}\ar@{.>}[ddddll]\\\\&&{0}\ar[dd]\ar[rr]&&{P_{k}}\ar[dd]^{(id,0)}\ar[rr]&&{0}\\\\{0}\ar[rr]&& {\oplus_{j\rightarrow k} P_{j}} \ar[rr]^{((\alpha_{j})_{j},0)}&&{\underline{P_{k}\oplus (\oplus_{j\rightarrow k} P_{j})}}\ar[rr]&&{0} }
\end{equation}
Similarly we can see it for the reverse composition and therefore (2.1) is isomorphic to the stalk complex $P_{k}$ in degree zero. 
\end{enumerate}
Hence, the complex is tilting iff we have $Hom_{K(P(Q))}(T^{k}_{k},T^{k}_{s}[-1])=0$, $\forall s$.
\end{proof}

\begin{definition}
Given a quiver with potential $(Q,S)$, define $\delta(Q,S)$ as the set of vertices for which the complex above is tilting over $J(Q,S)$, i.e.,
\begin{center}
$\delta(Q,S) = \left\{k\in Q_{0}: Hom_{K(P(Q))}(T^{k}_{k},T^{k}_{s}[-1])=0,\  \forall s \right\}$.
\end{center}
If $\delta(Q,S) \neq \emptyset$, then we say that $(Q,S)$ is \textbf{locally dualisable} in $\delta(Q,S)$. Furthermore, if $\delta(Q,S) = Q_{0}$ then we say that $(Q,S)$ is \textbf{globally dualisable}.\\
\end{definition}

\begin{remark}
Note that to check whether the complex is tilting we just need to check that there is no element $f$ in the path algebra such that 
\begin{equation}
\xymatrix{{\oplus_{j \rightarrow k} P_{j}}\ar[d]\ar[r]^{(\alpha_{j})_{j}}& {P_{k}}\ar[d]^{f}\\  {0}\ar[r] & {P_{s}}}
\end{equation}
commutes. The existence of such an $f$ implies that the set of relations must contain the set $\left\{f\alpha_{j} : j \to k \right\}$, which is easy to see once we differentiate the potential in order to the arrows. 
\end{remark}

The above remark allows us, given a potential $S$ for $Q$, to determine $\delta(Q,S)$.\newline From now on we'll drop the superscript on $T$ whenever the vertex with respect to which we're considering the tilting complex is fixed.\\
\begin{definition}
The \textbf{Seiberg dual} algebra of a quiver $Q$ with potential $S$ (or of its jacobian algebra) at the vertex $k\in \delta(Q,S)$ is the endomorphisms algebra of $T^{k}$ as defined above.
\end{definition}
Rickard's theorem then tells that seiberg dual algebras have derived equivalent categories of modules.
\end{section}
\begin{section}{Seiberg duality for good potentials}
Let's consider the following class of potentials:\\
\begin{definition}
A potential is said to be a \textbf{good potential} if each arrow appears at least twice and no subpath of length two appears twice.
\end{definition}
Note that, in particular, a quiver with a good potential has the property that every arrow is contained in at least two distinct cycles. \\
\begin{proposition}
A quiver with good potential is globally dualisable.
\end{proposition}

\begin{proof}
This is an immediate consequence of the definition of good potential since all the relations we get from these kind of potentials are of the form $\partial S/\partial a = \sum_{i=1}^{d}\lambda_{i}v_{i} = 0$ where $\lambda_{i}\in \mathbb{K}$, $d\geq2$ and the $v_{i}$'s are paths starting with different arrows thanks to the requirement that no subpath of length two should be shared between two terms of the potential. Thus, there cannot occur any relations of the type $u\alpha_{j} = 0$ and therefore $\delta(Q,S)=Q_{0}$.
\end{proof}

Let $(Q,S)$ be a quiver with good potential. We want to give a presentation of its Seiberg dual algebra at a fixed vertex $k$. We will see that this algebra is in fact the jacobian algebra of a quiver with potential. We will call this quiver the \textbf{Seiberg dual quiver}.

First we should compute the quiver. It has the same number of vertices as the initial quiver (since we will have that number of indecomposable projectives in $End_{D^{b}(Q)}(T)$ corresponds to the number of direct summands of $T$) and, for each irreducible homomorphism between the $T_{i}$'s, draw an arrow between the correspondent vertices. As we'll see in the next theorem, those irreducible homomorphisms are of three types (this terminology, used for simplicity of language, is inspired by [10]):

\begin{itemize}
\item arrows of the form $a$, where $a$ is also an arrow in $Q$, will be called \textbf{internal arrows}
\item arrows of the form $\alpha^{*}$ will be called \textbf{dual arrows};
\item arrows of the form $[\beta\alpha]$ will be called \textbf{mesonic arrows};
\end{itemize}
The theorem below shows that this choice of notation is an adequate one since the procedure to get the Seiberg dual quiver is the same as the one that allow us to mutate the initial quiver. Also as in mutations we will do this in two essential steps: obtain a quiver $\tilde{Q}$ that may contain more arrows than the irreducible homomorphisms and then, looking at relations, eliminate the appropriate arrows that do not correspond to irreducible ones (those will be the arrows lying in 2-cycles).\newline It turns out that relations on the Seiberg dual quiver can also be encoded in a potential (see Proposition 3.4) and it will be determined as follows:

\begin{enumerate}
 \item Determine $\tilde{S} := [S] + \sum_{i=1}^n [\beta_{i}\alpha_{j}]\alpha_{j}^{*}\beta_{i}^{*}$ (eventually containing some arrows representing non-irreducible homomorphisms);
 \item For every arrow $a$ in a two cycle $ab$, take the relation $\partial \tilde{S}/\partial a = 0$ and substitute $b$ in $\tilde{S}$ using this equality (and thus eliminate $b$ from the quiver, since $b$ is not irreducible as it can be written as a compostion of arrows). Call $\bar{S}$ to the potential thus obtained.\\
\end{enumerate}

\begin{remark} Again, for language simplicity, arrows appearing in two cycles will be called \textbf{massive arrows} and the process described on item 2 of the algorithm above will be called \textbf{integration over massive arrows}.
\end{remark}
Let us start by comparing the mutated quiver and the Seiberg dual quiver (no relations on them, yet).\\

\begin{theorem}
Let $Q$ be a quiver with a good potential $S$. The underlying quiver of the mutation of $Q$ coincides with the underlying quiver of the Seiberg dual of $Q$.
\end{theorem}
\begin{proof}
\begin{enumerate}
	\item First we prove that Seiberg duality at $k$ inverts incoming arrows to $k$. The complex $T_{k}$ has in degree zero one copy of $P_{j}$ for every arrow from $j$ to $k$, therefore for each such arrow you get one projection map from the direct sum to $P_{i}$ and therefore an irreducible homomorphism from $T_{k}$ to $T_{j}$, hence getting an arrow from $k$ to $j$ in the dual quiver. For each arrow $\alpha_{j}$ from $j$ to $k$, denote the correspondent homomorphism from $T_{k}$ to $T_{j}$ by $\alpha_{j}^{*}$. There are no more irreducible homomorphims: any other homomorphism factors through some factor of the direct sum first. 
	\item Now we prove that Seiberg duality at $k$ inverts outgoing arrows from $k$. This requires the commutativity of a diagram like the following:
\begin{equation}\nonumber
\xymatrix{{0}\ar[r]&{P_{i}}\ar[r]\ar[d]^{f}&{0}\ar[d]\\{0}\ar[r]&{\oplus_{j\rightarrow k} P_{j}}\ar[r]^{(\alpha_{j})_{j}}&{P_{k}}\ar[r]&{0}}
\end{equation}
The commutativity of the diagram requires that $(\alpha_{j})_{j} o f = 0$ and so we have to check the relations on the quiver to obtain such a condition. Fix an arrow $\beta$ from $i$ to $k$ and take the (cyclic) derivative of the potential in order to $beta$. Since $S$ is a good potential, $\partial S/\partial\beta = \sum_{t=1}^{d}\lambda_{t}v_{t}$ where the $v_{t}$'s  are paths from $i$ to $k$ (since $\beta v_{t}$ is a cycle for all $t$) and $d\geq2$. To give a homomorphism from $P_{i}$ to $\oplus_{j\rightarrow k} P_{j}$ we just need to give a homomorphism from $P_{i}$ to each $P_{j}$ by the universal property of the direct sum. Call $j_{t}$ to the index of the projective factor in $\oplus_{j\rightarrow k} P_{j}$ such that $\alpha_{j_{t}}$ is on the path $v_{t}$. Observe that $v_{t} = \alpha_{j_{t}}\tilde{v_{t}}$, where $\tilde{v_{t}}$ is a path from $i$ to $j_{t}$ as in the picture.
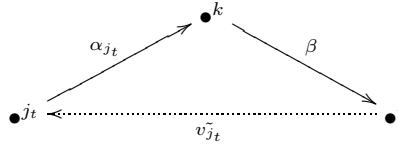
\begin{equation}\nonumber
\xymatrix{&&{\bullet}^{k}\ar[rrd]^{\beta}\\{\bullet}^{j_{t}}\ar[urr]^{\alpha_{j_{t}}}&&&&{\bullet}^{i}\ar@{.>}[llll]^{\tilde{v_{j_{t}}}} }
\end{equation}
Set the homomorphism from $P_{i}$ to each $P_{j}$ as follows:
\begin{itemize}
	\item zero if $j \neq j_{t}$ for some $t$;
	\item $\lambda_{t}\tilde{v_{t}}$ if $j = j_{t}$ for some $t$;
\end{itemize}
and set $\beta^{*}$ to be the homomorphism induced by this set of homomorphisms to the direct sum and therefore to the complex $T_{k}$. Clearly this map makes the diagram above commute. Now we need to prove that this is irreducible. If not, then it factors through other $T_{r}$ with the homomorphism from $T_{i}$ to $T_{r}$ being irreducible and therefore coming from an arrow $\gamma: i\rightarrow r$ ($r \neq k$ since the quiver has no two cycles). But the existence of such a factorization would imply that all terms $\beta v_{t}$ in the potential share a subpath of length two $\gamma\beta$ which is a contradiction since, by assumption, the potential is a good one and $d\geq2$. Hence $\beta^{*}$ is irreducible. By construction, these homomorphisms are the only irreducible ones.
\item For each path of length two in the initial quiver of the form \xymatrix{{\bullet}^{j}\ar[r]^{\alpha_{j}}&{\bullet}^{k}\ar[r]^{\beta_{i}}&{\bullet}^{i}}\  we get a homomorphism $\beta_{i}\alpha_{j}$ from $T_{j}$ to $T_{i}$. It will be irreducible iff there isn't a homomorphism in the opposite direction. In fact this follows from the fact that $S$ is a good potential and therefore every arrow appears in S. Thus, if $a$ is the arrow going in the opposite direction, $\partial\tilde{S}/\partial a$ gives an explicit factorization of the mesonic arrow. On the other hand, if it isn't contained in a 2-cycle, then it is irreducible, since it could only factor through the stalk complex of $P_{k}$ which is not, however, projective in $End_{D^{b}(Q)}(T)$. Denote this homomorphism by $[\beta_{i}\alpha_{j}]$.

\item Finally, if none of the previous cases apply, then homomorphisms between $T_{j}$ and $T_{i}$ are just arrows from $j$ to $i$. Again, these homomorphisms are irreducible iff they are not contained in a 2-cycle and a similar argument to the one above applies to this case.
\end{enumerate}
Let $\tilde{Q}$ be the quiver obtained by taking all the homomorphisms considered in the cases above, even if they are not irreducible. Determining this quiver $\tilde{Q}$ is, therefore, clearly the same procedure via mutations or via Seiberg duality. Now, since both mutation and Seiberg duality require the elimination of 2-cycles after this step (in the later case to get only the irreducible homomorphisms), the quiver obtained by mutation at $k$ and the quiver obtained by Seiberg duality on $k$ are the same.
\end{proof}

At this point, we shall prove that the algorithm above allows us to obtain the Seiberg dual potential of a quiver with potential $(Q,S)$ on a fixed vertex $k$.\\

\begin{proposition} The algorithm described above computes a potential for the Seiberg dual quiver such that its jacobian algebra is $End_{D^{b}(Q)}(T)$.
\end{proposition}
\begin{proof}
Let the homomorphisms represented by dual arrows of outgoing arrows be as it is described on the proof of (3.4). We will first prove that the relations induced by the potential $\tilde{S}$ obtained through the algorithm above are satisfied by $End_{D^{b}(Q)}(T)$. Let $\tau(i,j)$ be the coefficient of $[\beta_{i}\alpha_{j}]$ in $[S]$.\newline
\begin{itemize}
\item Relations coming from differentiating on $\beta_{i}^{*}$ (dual of an outgoing arrow):

\begin{equation}\nonumber
\partial \tilde{S} / \partial \beta_{i}^{*} = \sum_{j=1}^n [\beta_{i}\alpha_{j}]\alpha_{j}^{*} = \beta_{i}(\alpha_{j})_{j} = 0,
\end{equation}
since the map in question is homotopic to zero in the complex category.

\item Relations coming from differentiating on $\alpha_{j}^{*}$ (dual of an incoming arrow):

\begin{equation}\nonumber
\partial \tilde{S}/\partial \alpha_{j}^{*} = \sum_{i} \beta_{i}^{*}[\beta_{i}\alpha_{j}] = (\sum_{i} \beta_{i}^{*}\beta_{i})\alpha_{j}
\end{equation}
Let's check that $\sum_{i} \beta_{i}^{*}\beta_{i} = 0$. In fact, let's compute the m-th entry of this vector. For that we look to the appearences of $\alpha_{m}$ in $[S]$. So, if we have in $[S]$ some subexpression of the form
\begin{equation}\nonumber
\sum_{t=1}^{d}\tau(i_{t},m)[\beta_{i_{t}}\alpha_{m}]\tilde{v}_{i_{t}} 
\end{equation}
then we have the m-th entry of $\sum_{i} \beta_{i}^{*}\beta_{i}$ given by
\begin{equation}\nonumber
\sum_{t=1}^{d}\tau(i_{t},m)\tilde{v}_{i_{t}}\beta_{i_{t}}
\end{equation}
which is zero since $\partial S/\partial \alpha_{m} = 0$.

\item Relations coming from differentiating on $a$:

\begin{equation}\nonumber
\partial \tilde{S}/\partial a = \partial [S]/ \partial a = 0,
\end{equation}
since this is essentially the same as $\partial S/\partial a$ (eventually with some extra square brackets).

\item Relations coming from differentiating on $[\beta_{i}\alpha_{j}]$ (mesonic arrow):

We just need to check that 
\begin{equation}\nonumber
\partial [S]/\partial [\beta_{i}\alpha_{j}] = \alpha_{j}^{*}\beta_{i}^{*}
\end{equation}
but this follows by definition of $\alpha_{j}^{*}$ and $\beta_{i}^{*}$ as homomorphisms (see proof of (3.3))
\end{itemize}
Observe now that integration over massive arrows does not change the relations induced by the potential since the expressions obtained by differentiating in order to a massive arrows are zero in the jacobian algebra, according to the proof above.\newline
The last thing we need to check is that this potential $\bar{S}$ gives all the relations. If not, suppose first that the potential of $\bar{Q}$ is of the form $\bar{S} + W$. Then, for any arrow $a$ in W 
\begin{equation}\nonumber
\partial \bar{S}/\partial a + \partial W/\partial a = 0
\end{equation}
implies that $\partial W/\partial a = 0$ and hence $W = 0$. Suppose now that there is one non-zero relation $r$ such that $r$ is not of the form $\partial \bar{S}/\partial a$ for all $a$ in the quiver. This relation is a linear combination of homomorphisms such that each term of the linear combination is a map from some fixed $T_{j}$ to some fixed $T_{i}$. 
If this relation does not involve dual arrows, then these homomorphisms can be expressed as linear combinations of elements of the path algebra from $j$ to $i$ and therefore this is a relation iff there is some internal arrow $a$ such that $\partial S/\partial a$ is equal to $r$ up to square brackets. Thus we get a contradiction and therefore $r$ has to involve dual arrows. However, the construction of dual arrows as homomorphisms makes it easy to see that all possible relations involving them are contemplated in the cases above and thus proving that in fact all the relations are encoded in the potential $\bar{S}$.
\end{proof}

\begin{definition}
If one massive arrow appears in two different 2-cycles of $\tilde{S}$, that is, we get an expression of the form:
\begin{equation}\nonumber
\tilde{S} = \sum_{i=1}^{d}\lambda_{i}ab_{i} + \sum_{j=1}^{l}au_{j} + W
\end{equation}
where, $\lambda_{i}\in \mathbb{K}$; $a$ and $b_{i}$'s are arrows; $d \geq 2$; the $u_{i}$'s are paths of length $\geq 2$ and $a$ doesn't appear in W, then we say that the $b_{i}$'s are \textbf{related arrows}. 
\end{definition}

Given $Q$ a quiver with good potential $S$, suppose that $\tilde{S}$ can be written as follows:
\begin{equation}
\tilde{S} = \sum_{i=1}^{N} (\lambda_{i}a_{i}b_{i} + \sum_{j}\sigma_{i,j}a_{i}u_{i,j} + 	b_{i}v_{i}) + W 
\end{equation}
where $\sigma_{i,j}$, $\lambda_{i}$, $\mu_{i} \in \mathbb{K}$, the $a_{i}b_{i}$'s are 2-cycles (i.e., the $a_{i}$'s and the $b_{i}$'s are massive arrows), the $b_{i}$'s are mesonic (thus the coefficient of $b_{i}v_{i}$ is 1), W doesn't have any term involving massive arrows and $i\neq j$ implies  $a_{i} \neq a_{j}$  (that is, no related arrows occur). Note that $b_{i} \neq b_{j}$ because of the fact that $S$, being good, doesn't have repeated subpaths of length two.\\
\begin{theorem}
Let $Q$ be a quiver with a good potential $S$. If $k$ is a vertex such that no related arrows occur in the mutation, there is a right equivalence $\phi$ from $(\tilde{Q},\tilde{S})$ to $(\tilde{Q}, S' + \bar{S})$, where $S'$ is trivial and $\bar{S}$ is obtained by Seiberg duality.
\end{theorem}
\begin{proof}
Since there are no related arrows, let's assume that $\tilde{S}$ is of the form (3.1). Take the automorphisms given by:
\begin{equation}\nonumber
\begin{array}{rcl} \phi_{i}: \mathbb{K}\tilde{Q} & \to & \mathbb{K}\tilde{Q}  \\ a_{i}  & \mapsto & a_{i} - \frac{1}{\lambda_{i}} v_{i} \\b_{i} & \mapsto & b_{i} - \frac{1}{\lambda_{i}}\sum_{j}\sigma_{i,j}u_{i,j}\\ z & \mapsto & z \ \ \ \   $if$\ \ z \neq a_{i}$,$b_{i}$,\ $ z\in Q_{1} \end{array}
\end{equation}
Computing $\phi(\tilde{S})$, being $\phi$ the composition of all $\phi_{i}$'s, we get:
\begin{equation}\nonumber
\phi(\tilde{S}) = \sum_{i=1}^{N} (\lambda_{i}a_{i}b_{i} - \frac{1}{\lambda_{i}}\sum_{j}\sigma_{i,j}u_{i,j}v_{i}) + W 
\end{equation}
which reduced part is exactly 
\begin{equation}\nonumber
\sum_{i=1}^{N} (- \frac{1}{\lambda_{i}}\sum_{j}\sigma_{i,j}u_{i,j}v_{i}) + W
\end{equation}

Now, if we do integration over massive arrows in (3.1), taking in account that:
\begin{equation}\nonumber
\partial\tilde{S}/\partial a_{i} = \lambda_{i}b_{i} + \sum_{j}\sigma_{i,j}u_{i,j} \ \ \ \ \ \ 
\partial\tilde{S}/\partial b_{i} = \lambda_{i}a_{i} + v_{i} 
\end{equation}
and using the relations $\partial\tilde{S}/\partial a_{i} = 0$ and $\partial\tilde{S}/\partial b_{i} = 0$ in $\tilde{S}$ we get:
\begin{equation}\nonumber
\sum_{i=1}^{N}(- \frac{1}{\lambda_{i}}\sum_{j}\sigma_{i,j}u_{i,j}v_{i}) + W
\end{equation} 
which is the same as $\phi(\tilde{S})_{red}$.
\end{proof}

\begin{corollary}
If $Q$ is a quiver with a good potential $S$ and if $k$ is a vertex such that no related arrows arise in the mutation procedure, then mutation at $k$ and Seiberg duality at $k$ are isomorphic. In particular mutation of good potentials give derived equivalent algebras.
\end{corollary}
\begin{proof}
We have that $J(Q_{triv}\oplus Q_{red}, S_{triv} + S_{red}) \cong J(Q_{red}, S_{red})$. Conjugating this fact with remark 1.5 and theorems 3.4 and 3.7, we get the result.
\end{proof}

\end{section}
\begin{section}{An example}
Given a Del Pezzo surface $S$, we can realise its derived category of coherent sheaves as a path algebra with relations. This can be done using strongly exceptional collections. For the purpose of what follows, let us recall some results and definitions.\\

\begin{definition}
An \textbf{exceptional collection} on a projective surface is a collection of coherent sheaves $\left\{E_{1}, ..., E_{n}\right\}$ such that:
\begin{itemize}
	\item Ext$^{k}$(E$_{i}$,E$_{i}$)=0, $\forall k > 0$ and Hom(E$_{i}$,E$_{i}$)= $\mathbb{K}$
	\item Ext$^{k}$(E$_{i}$,E$_{j}$)=0, $\forall  1 \leq j < i \leq n$, $\forall k > 0$
	\item The stalk complexes of these sheaves generate D$^{b}$(Coh(X)) as a triangulated category.
	\end{itemize}
It is \textbf{strongly exceptional} if, furthermore, Ext$^{k}$(E$_{i}$,E$_{j}$)=0, $\forall  1 \leq i,j \leq n$, $\forall k > 0$\\ \\
\end{definition}

\begin{theorem}
If $S$ is a Del Pezzo Surface, we have a strongly exceptional collections of sheaves given by:
\begin{itemize}
	\item $\left\{O, O(1), O(2)\right\}$ if S = $\mathbb{P}^{2}$
	\item $\left\{O, O(1,0), O(0,1), O(1,1)\right\}$ if S = $\mathbb{P}^{1} \times \mathbb{P}^{1}$
	\item $\left\{O, O(E_{1}), ..., O(E_{r}), O(1), O(2)\right\}$ if S is $dP_{r}$ with $r\leq 8$, where each $E_{i}$ is an exceptional curve of the blow up and $dP_{r}$ is the Del Pezzo obtained by blowing up $1\leq r \leq 8$ points in $\mathbb{P}^2$.\\
\end{itemize}
\end{theorem}

\begin{definition}
Let $X$ be a nonsingular projective variety. A coherent sheaf $T$ is said to be \textbf{tilting} if:
\begin{itemize}
	\item $Ext^{k}(T,T)=0,  \forall k > 0$
	\item T generates $D^b(Coh(X))$ as triangulated category
	\item B = End(T) has finite global dimension\\
\end{itemize}
\end{definition}

\begin{theorem}
Let $X$ be a nonsingular projective variety, $T$ a coherent sheaf on $X$ and $B = End(T)$. Then the following are equivalent:
\begin{enumerate}
	\item T is tilting;
	\item There is an equivalence $\Phi: D^b(Coh(X)) \rightarrow D^b(mod(B))$ of triangulated categories with $\Phi(T) = B$, where mod(B) is the category of finitely generated right modules over B.
\end{enumerate}
\end{theorem}

Now, if we take the direct some of a strongly exceptional collection over $S$, we get a tilting sheaf and, therefore, a derived equivalence between $Coh(S)$ and $\mathbb{K}Q/I$ for some quiver $Q$ and some ideal of relations $I$. These are determined looking at the irreducible homomorphisms between the sheaves in the collection and taking relations between those homomorphisms.\\ 

\begin{example}
$dP1$ with exceptional collection $\left\{O, O(E_{1}), O(1), O(2)\right\}$
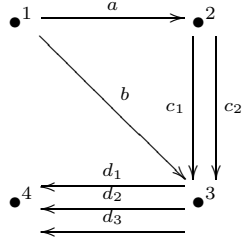
\begin{equation}\nonumber
\xymatrix{
{\bullet}^1\ar[rr]^{a}\ar[ddrr]^{b}&&{\bullet}^2\ar@<-1ex>[dd]_{c_{1}}\ar@<1ex>[dd]^{c_{2}}\\ \\ 
{\bullet}^4 &&{\bullet}^3\ar@<-1ex>[ll]_{d_{1}}\ar@<1ex>[ll]_{d_{2}} \ar@<3ex>[ll]_{d_{3}} }
\end{equation}
with relations:
\begin{equation}\nonumber
d_{3}c_{1} = d_{1}c_{2}
\end{equation}
\begin{equation}\nonumber
d_{2}c_{1}a = d_{1}b
\end{equation}
\begin{equation}\nonumber
d_{3}b = d_{2}c_{2}a
\end{equation}
\end{example}

This example, however, doesn't fit in our setting of quivers with potentials. In fact what we ought to consider is not $S$ itself but $X = \omega_S$ - the total space of the canonical bundle of $S$ - instead. This is a local Calabi-yau three-fold and if we let $\pi: X \rightarrow S$ be the natural projection, we get $\tilde{B}= End_{X}(\oplus_{i}\pi^{*}E_{i})$ is derived equivalent to $Coh(X)$, where $(E_{i})_{i}$ is an exceptional collection over S. This algebra $\tilde{B}$ can also be seen as a path algebra of a quiver which can be obtained from the correspondent Del Pezzo quiver adding one arrow for each relation in the opposite direction of the composition of arrows in that relation. These will be quivers with potentials.\\
\\
\begin{example}
The completed quiver for $dP1$ with exceptional collection 
\newline$\left\{O, O(E_{1}), O(1), O(2)\right\}$ is:
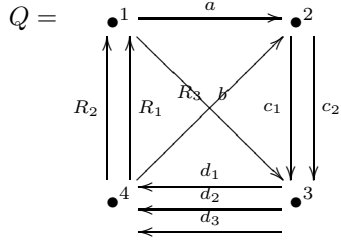
\begin{equation}\nonumber
Q = \xymatrix{
{\bullet}^1\ar[rr]^{a}\ar[ddrr]^{b}&&{\bullet}^2\ar@<-1ex>[dd]_{c_{1}}\ar@<1ex>[dd]^{c_{2}}\\ \\ 
{\bullet}^4\ar[uurr]^{R_{3}}\ar@<-1ex>[uu]_{R_{1}}\ar@<1ex>[uu]^{R_{2}} &&{\bullet}^3\ar@<-1ex>[ll]_{d_{1}}\ar@<1ex>[ll]_{d_{2}} \ar@<3ex>[ll]_{d_{3}} }
\end{equation}
with potential:
\begin{equation}\nonumber
S = R_{3}(d_{3}c_{1} - d_{1}c_{2}) + R_{1}(d_{1}b - d_{2}c_{1}a) + R_{2}(d_{2}c_{2}a - d_{3}b)
\end{equation}

Note that this is a good potential. Therefore, using our results, mutations on any vertex of this quiver will give us derived equivalent path algebras. Since the one above is derived equivalent to $Coh(X)$, so will be $\mu_k(Q,S)$. Let's finish by presenting $\mu_1(Q,S)$.

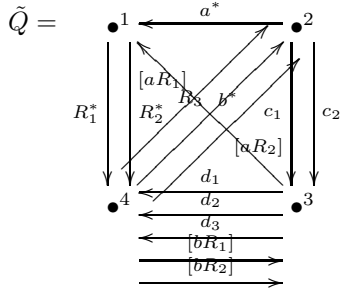
\begin{equation}\nonumber
\tilde{Q} = \xymatrix{
{\bullet}^1\ar@<-1ex>[dd]_{R_{1}^{*}}\ar@<1ex>[dd]^{R_{2}^{*}}&&{\bullet}^2\ar[ll]_{a^{*}}\ar@<-1ex>[dd]_{c_{1}}\ar@<1ex>[dd]^{c_{2}}\\ \\ 
{\bullet}^4\ar@<-2ex>[uurr]_{[aR_{2}]}\ar@<2ex>[uurr]^{[aR_{1}]}\ar[uurr]^{R_{3}}\ar@<-5ex>[rr]^{[bR_{1}]}\ar@<-7ex>[rr]^{[bR_{2}]} &&{\bullet}^3\ar[uull]_{b^{*}}\ar@<-1ex>[ll]_{d_{1}}\ar@<1ex>[ll]_{d_{2}} \ar@<3ex>[ll]_{d_{3}} }
\end{equation}

We take a cyclically equivalent potential since there are terms on it starting and ending at $1$. Then let's substitute paths of length two passing through $1$ by new arrows and add $\Delta_{1}$. 
\begin{equation}\nonumber
\tilde{S} = R_{3}d_{3}c_{1} - R_{3}d_{1}c_{2} - d_{2}c_{1}[aR_{1}] + d_{1}[bR_{1}] - d_{3}[bR_{2}] + d_{2}c_{2}[aR_{2}] 
\end{equation}
\begin{equation}\nonumber
+ [aR_{1}]R_{1}^{*}a^{*} + [aR_{2}]R_{2}^{*}a^{*} + [bR_{1}]R_{1}^{*}b^{*} + [bR_{2}]R_{2}^{*}b^{*}
\end{equation}
Clearly this potential is not reduced. Let's condsider the following right equivalence:
\begin{equation}\nonumber
\begin{array}{rcl} \phi: \mathbb{K}\tilde{Q} & \to & \mathbb{K}\tilde{Q} \\ d_{1}  & \mapsto & d_{1} - R_{1}^{*}b^{*} \\d_{3} & \mapsto & -d_{3} + R_{2}^{*}b^{*}\\ \left[bR_{1}\right] & \mapsto & \left[bR_{1}\right] + c_{2}R_{3} \\ \left[bR_{2}\right] & \mapsto & \left[bR_{2}\right] + c_{1}R_{3} \\ u & \mapsto & u \   \   $if$\  u \neq d_{1},d_{3},[bR_{1}],[bR_{2}]$,$\ \  u\in Q_{1} \end{array}.
\end{equation}
If we compute $\phi(\tilde{S})$, it is of the form $S' + \bar{S}$ and thus we can take the reduced part or integrate over massive arrows. In any case, as proved in Theorem 3.7, we get the same result which is:
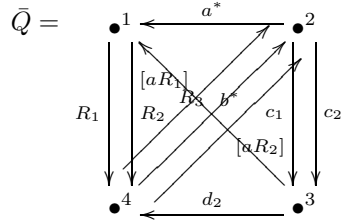
\begin{equation}\nonumber
\bar{Q} = \xymatrix{
{\bullet}^1\ar@<-1ex>[dd]_{R_{1}}\ar@<1ex>[dd]^{R_{2}}&&{\bullet}^2\ar[ll]_{a^{*}}\ar@<-1ex>[dd]_{c_{1}}\ar@<1ex>[dd]^{c_{2}}\\ \\ 
{\bullet}^4\ar@<-2ex>[uurr]_{[aR_{2}]}\ar@<2ex>[uurr]^{[aR_{1}]}\ar[uurr]^{R_{3}} &&{\bullet}^3\ar[uull]_{b^{*}}\ar@<1ex>[ll]_{d_{2}}  }
\end{equation}
with potential
\begin{equation}\nonumber
\bar{S} = c_{2}R_{3}R_{1}^{*}b^{*} + c_{1}R_{3}R_{2}^{*}b^{*} + d_{2}c_{2}[aR_{2}] - d_{2}c_{1}[aR_{1}] + [aR_{1}]R_{1}^{*}a^{*} + [aR_{2}]R_{2}^{*}a^{*}.
\end{equation}
Since this new jacobian algebra is derived equivalent to $J(Q,S)$, it is also derived equivalent to $Coh(X)$.
\end{example}

\end{section}

\end{document}